\documentclass{article}
\usepackage{amsmath,amssymb,hyperref}
\usepackage{graphicx}
\usepackage{subcaption}
\usepackage{geometry}

\newcommand{\matlab}{{\sc matlab}}
\newcommand{\imagunit}{\mathbf{i}}
\newcommand{\C}{{\mathbb{C}}}
\newcommand{\R}{{\mathbb{R}}}
\newcommand{\Mn}{{\mathbf{M}_{n}}}
\newcommand{\lmax}{\lambda_{\max}}

\title{An Experimental Comparison of Methods for Computing the Numerical Radius}

\author{Tim Mitchell\thanks{
Department of Computer Science, Queens College/CUNY,  
\href{mailto:tim.mitchell@qc.cuny.edu}{\texttt{tim.mitchell@qc.cuny.edu}}}
 \and Michael L. Overton\thanks{
 Courant Institute of Mathematical Sciences, New York University,
 \href{mailto:mo1@nyu.edu}{\texttt{mo1@nyu.edu}}}
 }
 \date{October 6, 2023}
 \begin{document}
 \maketitle
 \begin{abstract}
 We make an experimental comparison of methods for computing the numerical radius
 of an~\mbox{$n\times n$} complex matrix,
 based on two well-known characterizations, the first a nonconvex optimization problem in
 one real variable and the second a convex optimization problem in~\mbox{$n^{2}+1$} real variables. 
 We make comparisons with respect to both accuracy and computation time using
 publicly available software.
 \end{abstract}
 
 \small
 \centerline{AMS Subject Classifications: 15A60, 90C22}
 
\section{Introduction}
 Let $A\in\Mn$, the space of $n\times n$ complex matrices. The numerical radius $r(A)$ is defined
 as the maximum, in modulus, of the points in the field of values
 (numerical range), i.e.,
 \[
        r(A) = \max\big \{ |v^{*}Av|:  v\in \C^{n}, \|v\|_{2} = 1 \big\}.
\]
A well-known result (Johnson \cite{Joh78}; see also Kippenhahn \cite{Kip51})
is
\begin{equation}\label{rdef}
    r(A) = \max_{\theta \in [0,2\pi)}  h(\theta)
\end{equation}
where 
\[
        h(\theta ) = \lmax(H(\theta )) 
\]
and 
\[
        H(\theta) =\frac{1}{2}\left( e^{\imagunit \theta} A + e^{-\imagunit \theta} A^{*}\right),
\]
with $\lmax(\cdot)$ denoting largest eigenvalue of a Hermitian matrix,
 $\imagunit$ the imaginary unit and $^{*}$ denoting complex conjugate transpose.
Although the optimization problem \eqref{rdef} has only one real variable, it is nontrivial
to solve because $h$ is nonconvex.  Various algorithms have been introduced to solve~\eqref{rdef}, 
including the level-set method of Mengi and Overton~\cite{MenOve05}
(based on earlier work of He and Watson~\cite{HeW97} and of Boyd and Balakrishnan~\cite{BoyB90}), and
a cutting plane method of Uhlig~\cite{Uhl09} (based on earlier work of Johnson~\cite{Joh78}).
More recently, improvements to both of these approaches have been
given by Mitchell~\cite{Mit23}, along with
a state-of-the-art hybrid method that combines the virtues of both in order
to remain efficient in all cases.

A second well-known characterization of $r(A)$ (Mathias \cite{Mat93}; see also Ando \cite{And73})
is
 \begin{equation}\label{SDPchar} 
       r(A)=\min_{c\in \R, Z \in \Mn,Z=Z^{*}}\left\{c:  \begin{bmatrix} cI_n +Z&A\\A^*& cI_n -Z\end{bmatrix}\succeq 0\right \},
\end{equation}
where $M \succeq 0$ means a Hermitian matrix
$M$ is positive semidefinite. Thus, $r(A)$ can be computed by solving a convex optimization problem, known as a semidefinite program (SDP),
in $n^{2} + 1$ real variables. If $A$ is real, then $Z$ can be restricted to be real symmetric, so the 
optimization problem has \mbox{$n(n+1)/2 + 1$} variables. It is well known that
such convex optimization problems can be solved, up to any given accuracy, in polynomial time
using the Turing complexity model;
the traditional method guaranteeing such complexity is the ellipsoid method, but more
practical interior-point methods have also recently been designed with polynomial-time complexity
\cite{deKVal16}.

However, for all practical purposes,
interior-point methods for computing $r(A)$ via \eqref{SDPchar}
are very slow compared to more direct methods
based on \eqref{rdef}. In this note
we illustrate this with some experimental computations carried out using 
Mitchell's codes available from the supplemental materials published in~\cite{Mit23} as well as
the latest stable release (version 2.2, build 1148) of 
the well-known CVX software for ``disciplined convex programming'' \cite{cvx14}.
CVX uses a variety of primal-dual interior-point ``solvers'' for SDP; we use 
two of these, SDPT3 and SeDuMi.
CVX lacks the polynomial-time guarantees of the methods described by
\cite{deKVal16}, but  it is much faster than such methods, partly because
it uses floating point arithmetic, while any rigorous polynomial-time method would necessarily
have to be implemented using arbitrary precision integer or rational
arithmetic, as discussed at the end of
\cite{deKVal16}.
 
 \section{The experiments and the results} 
We experimented with the following methods using \matlab\ R2023a Update 4 and 
 a 2023 Mac Mini configured with a 12-core M2 Pro Apple Silicon CPU and 16GB RAM running macOS 13.5:

\begin{figure}[t]
\centering
\subcaptionbox{Real Matrices}{\includegraphics[scale=0.415,trim=0cm 0cm 1.0cm 0cm,clip]{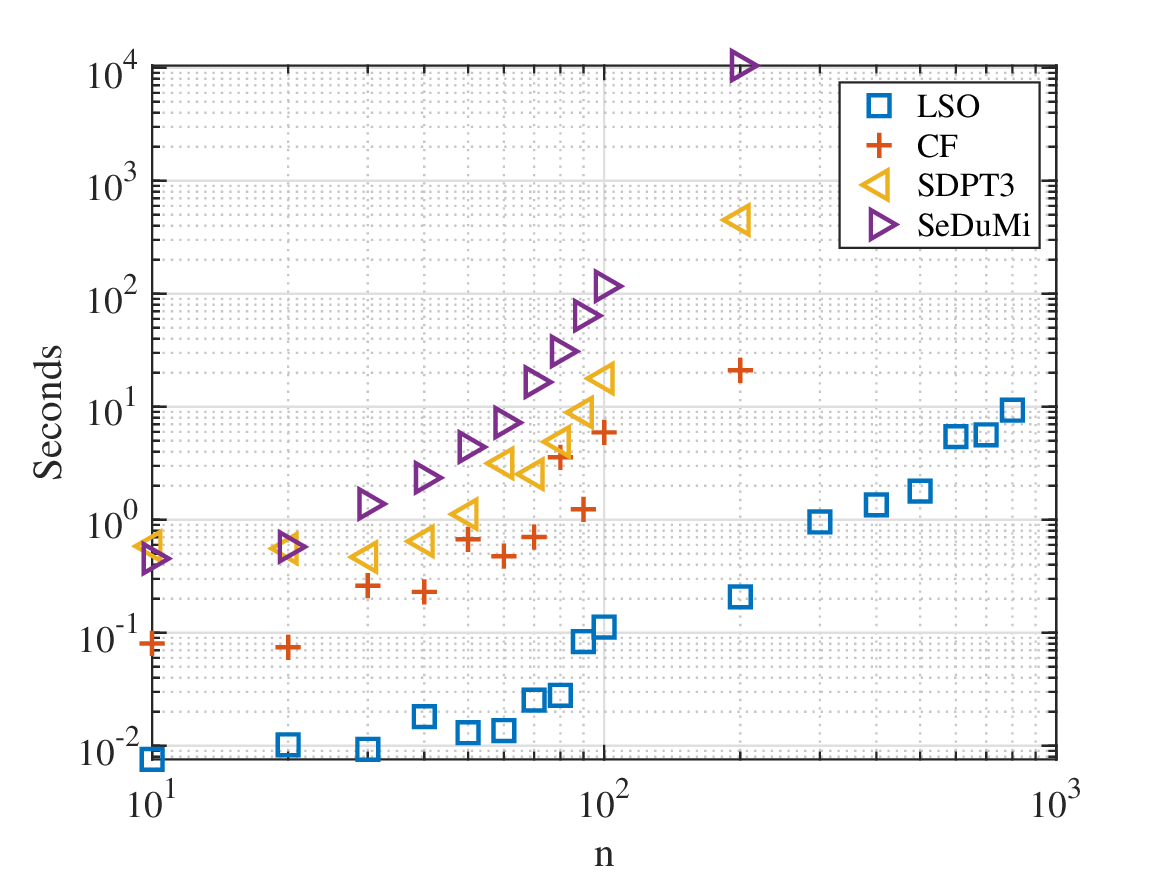}}%
\subcaptionbox{Complex Matrices}{\includegraphics[scale=0.415,trim=0.5cm 0cm 0cm 0cm,clip]{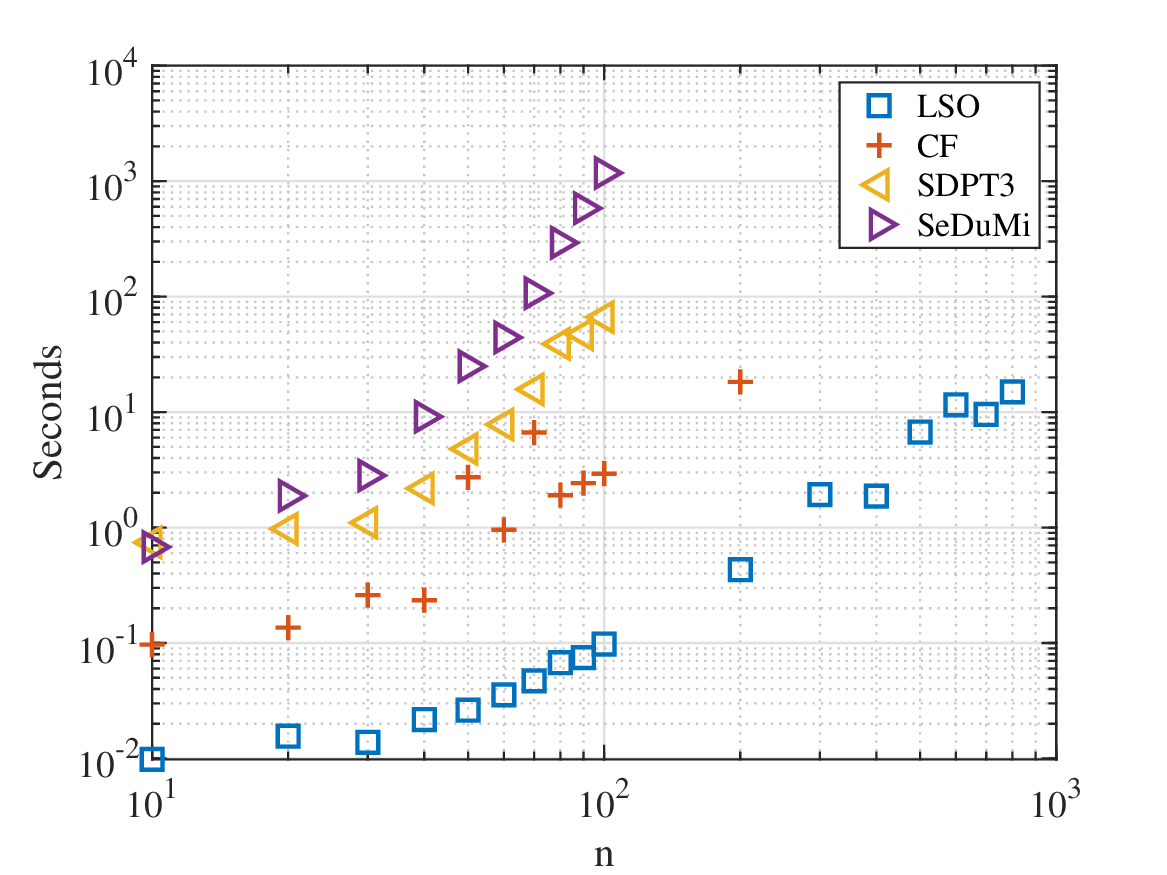}}%
\caption{The elapsed running times of the four methods tested are shown.
For problems that took at most 3 minutes to solve, the reported times are the average of five trials in order to account for
timing variability, while all other reported times are from a single run.
Due to their high cost, CF, SDPT3, and SeDuMi were only tested up to size~$n=200$, and
on the complex matrix with~$n=200$, both SDPT3 and SeDuMi crashed, and thus their times are only reported up to $n=100$ 
in the right plot.
}
\label{fig:times}
\end{figure} 
 
 \begin{itemize}
 \item LSO: solve \eqref{rdef} via the Level Set with Optimization Algorithm 3.1 in \cite{Mit23}\footnote{We do not include 
 Mitchell's improved
 cutting-plane method, Algorithm 5.1 in \cite{Mit23},
 or his hybrid method that combines Algorithms 3.1 and 5.1, as it is not necessary in our benchmark.
 Although Algorithm 5.1 can be many times faster than the Algorithm 3.1 on certain problems and many times slower on others,
 on the random matrices that we test here, Algorithm 3.1 and Algorithm 5.1 generally perform similarly.}
 \item CF: solve \eqref{rdef} using Chebfun via adaptive Chebyshev interpolation\footnote{Chebfun is a package for 
globally approximating functions to machine precision accuracy using adaptive Chebyshev interpolation techniques;
for more details, see~\cite{DriHT14}.
By method CF, we mean using Chebfun to accurately approximate $h$ on $[0,2\pi]$ 
followed by a call to the Chebfun {\tt max} routine to find the
 maximal value, as proposed in~\cite{Uhl09}.}
 \item SDPT3: solve \eqref{SDPchar} using CVX \cite{cvx14} with the SDPT3 solver
 \item SeDuMi:  solve \eqref{SDPchar} using CVX \cite{cvx14} with the SeDuMi solver 
 \end{itemize}

\begin{figure}
\centering
\subcaptionbox{Real Matrices}{%
	\includegraphics[scale=0.415,trim=1.9cm 0cm 3.7cm 0cm,clip]{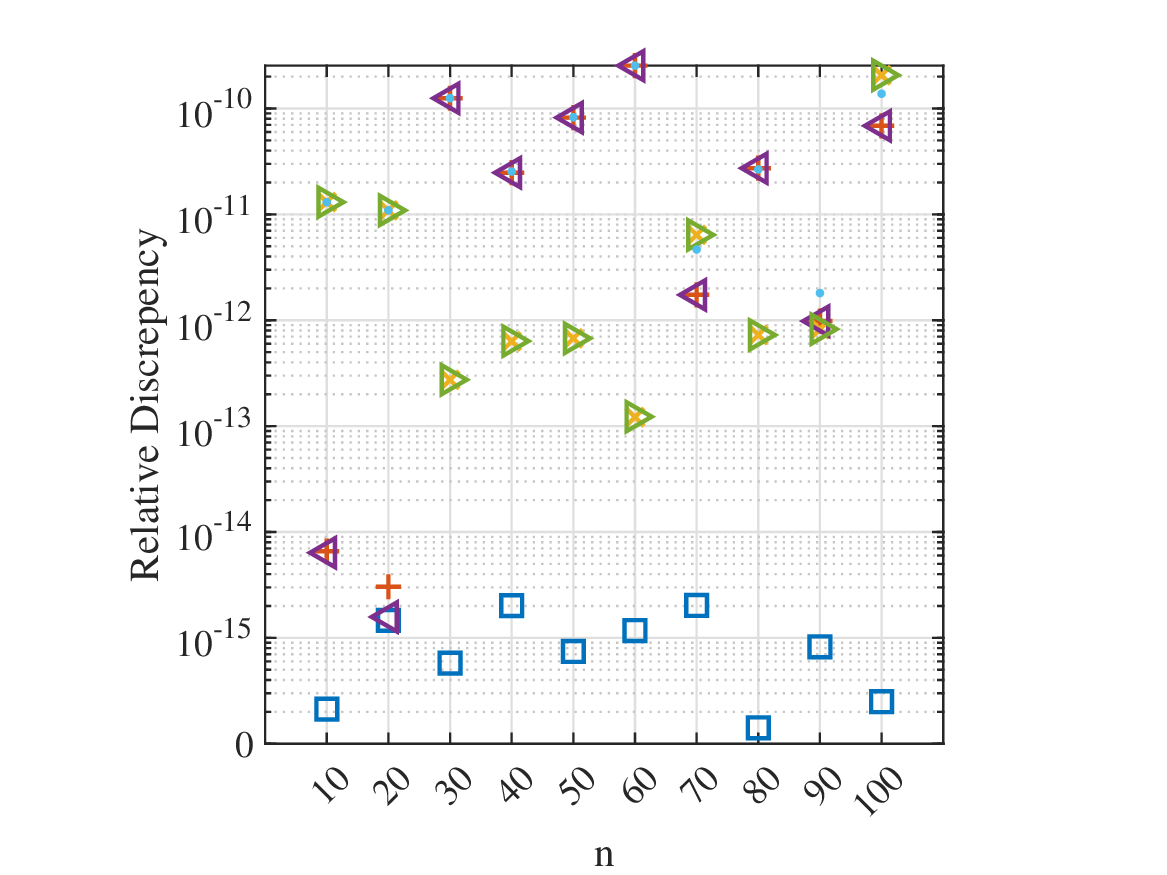}%
	\includegraphics[scale=0.453,trim=12.3cm 0cm 1.2cm 0cm,clip]{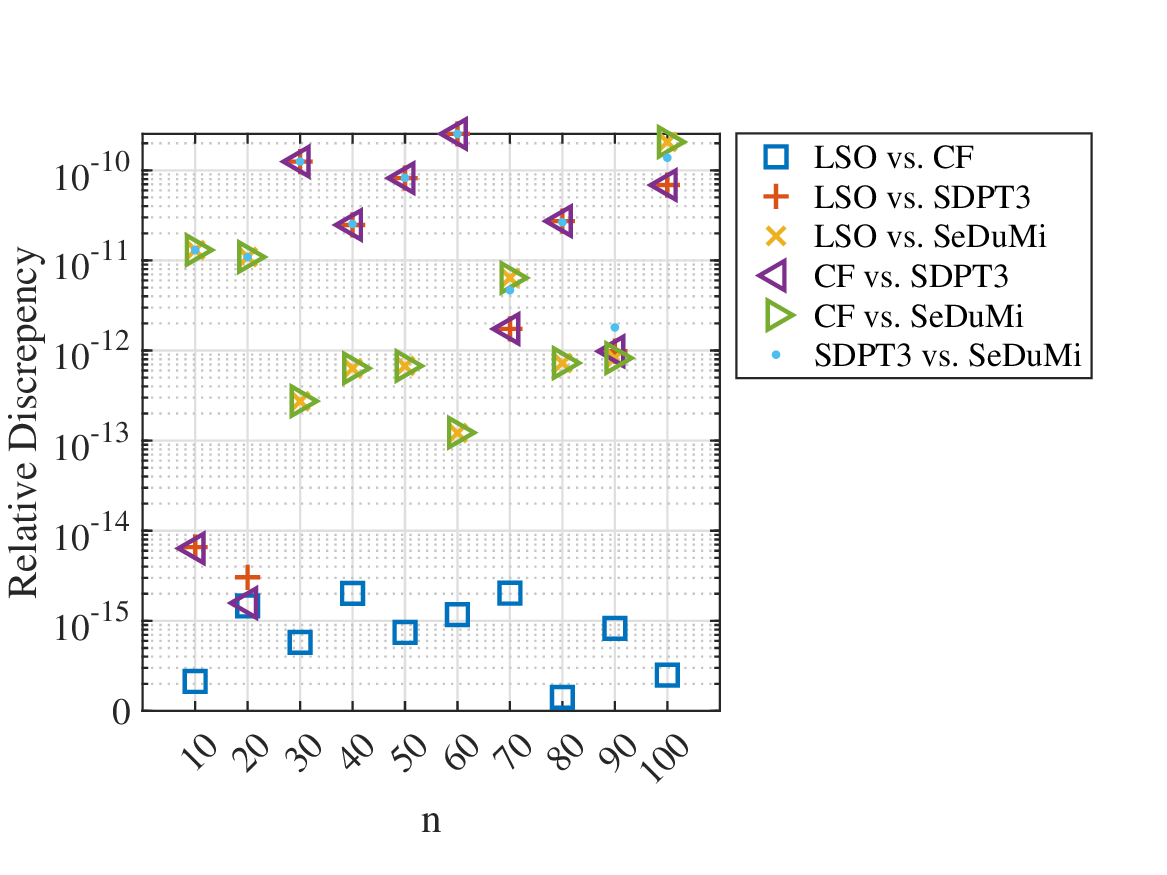}
}%
\subcaptionbox{Complex Matrices}{%
	\includegraphics[scale=0.415,trim=1.9cm 0cm 0cm 0cm,clip]{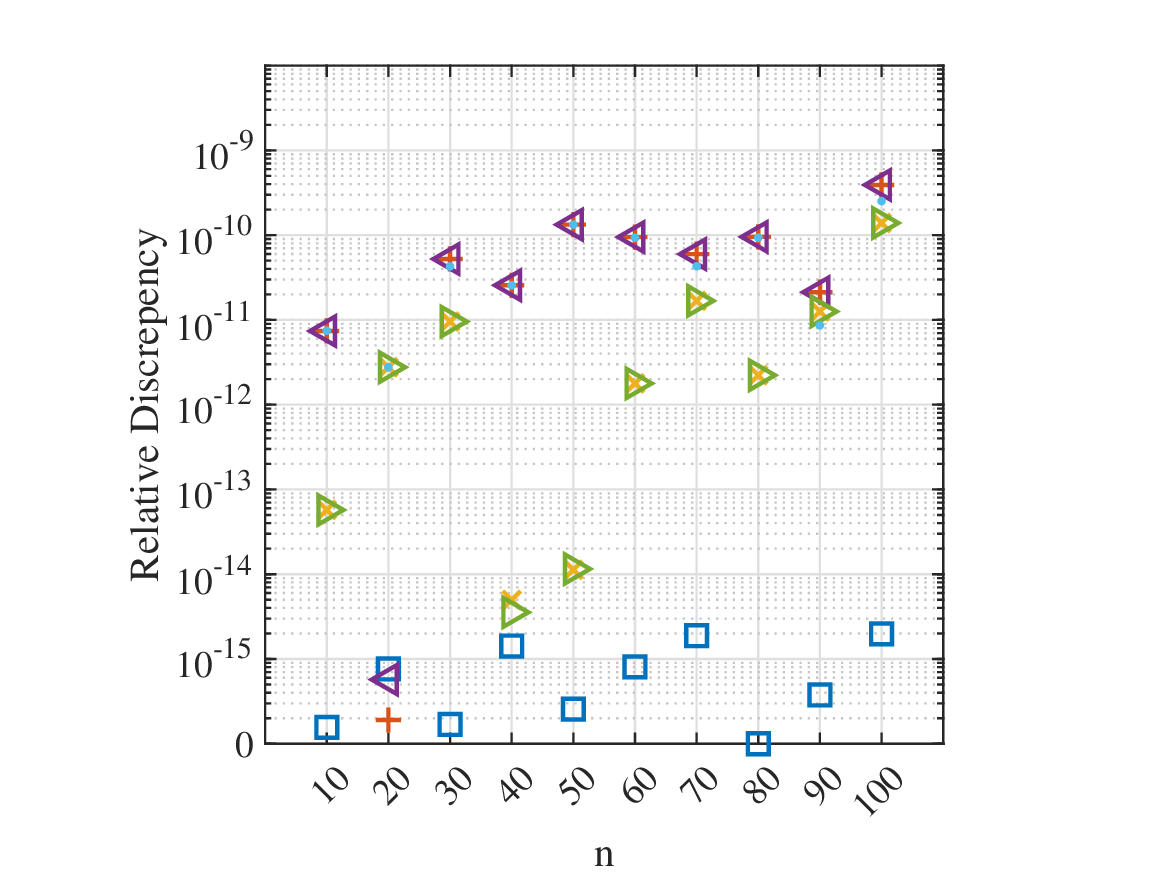}
}%
\caption{The plots give the pairwise relative discrepancies of the methods' computed numerical radius estimates, i.e., $|(r_{1}(A)-r_{2}(A))|/\max(r_{1}(A),r_{2}(A))$, where
$r_{1}(A)$ and $r_{2}(A)$ are the computed results for the two methods being compared
on matrix $A$.}
\label{fig:errors}
\end{figure}

We conducted our experiments using randomly 
generated real as well as complex matrices ranging in size from $n=10$ to $n=800$; 
the matrix coefficients (real parts and imaginary parts, as appropriate)
were generated using the standard normal distribution.  
Figure \ref{fig:times} shows the running times for our test matrices, 
with the left panel showing
 the results for real matrices and the right panel the results for complex
 matrices. We see immediately that
 LSO is far more efficient than the other three methods, and that among the other three, CF is the fastest, then SDPT3, and then SeDuMi.\footnote{We note that CVX uses
precompiled binary files for both SDPT3 and SeDuMi; 
if these were instead implemented in MATLAB, they'd presumably be signficantly slower.}  
 Meanwhile, Figure~\ref{fig:errors} shows, for each 
 pairing of the four methods, the relative discrepancies between the computed values 
 of the numerical radius $r$. We see that methods LSO and CF have relative
 discrepancies always near $10^{-15}$, indicating that these methods compute $r$ to
 about 14 correct digits, quite impressive considering that \matlab's IEEE double format
 numbers have only 53 bits of precision, or approximately 16 decimal digits. The discrepancies between the LSO method and the CVX methods, as well as between
 the CF method and the CVX methods, indicate that the CVX methods solving \eqref{SDPchar}
 return results which are accurate to only about~9~to~12 digits, although we note that CVX
 was instructed to compute the results with high precision accuracy.\footnote{By setting
 {\tt cvx\_precision high} in the calling code.}

  \section{Conclusion}
 
 Our experiments show that although the characterization of the numerical radius as the solution
 of a semidefinite program is a valuable theoretical tool, as a practical matter, Mitchell's
 Algorithm 3.1
 based on the nonconvex characterization is much faster. We note, however, that the SDP 
 formulation has the appealing property that it also allows the convenient solution of linearly
 parameterized optimization problems involving the numerical radius, 
 such as the one given in \cite[Sec.~2.1]{LewOve20}.

\bibliography{refs}
\bibliographystyle{alpha}
\end{document}